\documentclass{article}
\usepackage{amssymb,latexsym,amsmath,amsthm,mathrsfs}
\usepackage{graphicx}
\newcommand{\comment}[1]{}

\begin{document}
\title{An arithmetic theorem and its demonstration\footnote{Originally published as
{\em Theorema arithmeticum eiusque demonstratio},
Commentationes arithmeticae collectae
\textbf{2} (1849), 588--592.
E794 in the Enestr{\"o}m index.
Translated from the Latin by Jordan Bell,
Department of Mathematics, University of Toronto, Toronto, Ontario, Canada.
Email: jordan.bell@gmail.com}}
\author{Leonhard Euler}
\date{}
\maketitle

The Theorem which I have taken here to propose and demonstrate I had
previously communicated to colleagues,\footnote{Translator: Euler's
September 25, 1762 and November 9, 1762 letters to Goldbach.}
in which it was seen to be fully worthy of attention, especially
since its demonstration is hardly obvious and may perhaps have been
investigated by many in vain. I stated the theorem in the following
way:\footnote{Translator: Euler showed this in
his {\em Institutionum calculi integralis volumen secundum},
1769, E366, \S 1169, which is explained
in Ed Sandifer's March, 2005 {\em How Euler Did It}.
The {\em Opera omnia} also refers to
Euler's 
E540, {\em Nova methodus fractiones quascumque rationales in fractiones simplices resolvendi} (1775) and
E475, {\em Speculationes analyticae} (1774).}

{\em If however many distinct numbers $a,b,c,d$ etc. are given, and from
them fractions are formed whose common numerator is unity, and such that
the denominator of each is the product of all the differences of the number
and each of the remaining numbers, so that the fractions are
}
\[
\frac{1}{(a-b)(a-c)(a-d) \; \textrm{etc.}}, \;
\frac{1}{(b-a)(b-c)(b-d) \; \textrm{etc.}}, \;
\frac{1}{(c-a)(c-b)(c-d) \; \textrm{etc.}} \;
\textrm{etc.},
\]
{\em then the sum of all these fractions is always equal to $0$.}

Thus for example for the given numbers $2,5,7,8$, the four fractions
\[
\frac{1}{-3\cdot -5\cdot -6}, \quad \frac{1}{3\cdot -2\cdot -3},
\quad \frac{1}{5\cdot 2\cdot -1}, \quad \frac{1}{6\cdot 3\cdot 1}
\]
are to be thence formed, which reduce to these
\[
-\frac{1}{3\cdot 5\cdot 6}, \quad +\frac{1}{3\cdot 2\cdot 3},
\quad -\frac{1}{5\cdot 2\cdot 1},\quad +\frac{1}{6\cdot 3\cdot 1},
\]
and on the strength of the theorem
\[
-\frac{1}{90}+\frac{1}{18}-\frac{1}{10}+\frac{1}{18}=0.
\]
So as not to make complications with negative signs, these fractions can
be formed so that, arranging the given terms according to order
of magnitude, either increasing or decreasing, for each the differences
with it and all the other terms are taken, and taking these for the
denominators in fractions
with numerator unity, the
signs are taken alternately $+$ and $-$. 

For instance, if the given numbers are
\[
3, \; 8, \; 12, \; 15, \; 17, \; 18,
\]
from each the denominators are thus collected:
\[
\begin{array}{r|rl}
3&5\cdot 9\cdot 12\cdot 14\cdot 15&=113400\\
8&5\cdot 4\cdot 7\cdot 9\cdot 10&=12600\\
12&9\cdot 4\cdot 3\cdot 5\cdot 6&=3240\\
15&12\cdot 7\cdot 3\cdot 2\cdot 3&=1512\\
17&14\cdot 9\cdot 5\cdot 2\cdot 1&=1260\\
18&15\cdot 10\cdot 6\cdot 3\cdot 1&=2700
\end{array}
\]
and it will be
\[
\frac{1}{113400}-\frac{1}{12600}+\frac{1}{3240}-\frac{1}{1512}+\frac{1}{1260}
-\frac{1}{2700}=0
\]
or multiplying each by 36
\[
\frac{1}{3150}-\frac{1}{350}+\frac{1}{90}-\frac{1}{42}+\frac{1}{35}-\frac{1}{75}=0,
\]
and reducing these fractions to the single denominator 3150, it is
clear by itself that
\[
\frac{1-9+35-75+90-42}{3150}=0.
\]

Indeed, in the case in which only two numbers are given the theorem
needs no demonstration, since it is transparent that
\[
\frac{1}{a-b}+\frac{1}{b-a}=0;
\]
but even in the case of three numbers $a,b,c$ it is now more
subtle, for it is not immediately clear that
\[
\frac{1}{(a-b)(a-c)}+\frac{1}{(b-a)(b-c)}+\frac{1}{(c-a)(c-b)}=0;
\]
and for larger numbers, and even more so in general for any number whatsoever,
finding the truth for simpler cases offers scant help.

Indeed, I have extended this theorem more widely, and it can be
stated in the following way.

\begin{center}
{\Large A General Theorem\footnote{Translator: Euler in fact showed this in
his {\em Institutionum calculi integralis volumen secundum},
1769, E366, \S 1169.}}
\end{center}

{\em If however many distinct numbers $a,b,c,d,e,f$ etc. are given,
whose number $=m$, and the following products are formed from the differences
of each with the others}
\begin{eqnarray*}
(a-b)(a-c)(a-d)(a-e)(a-f)\;\textrm{etc.}&=&A,\\
(b-a)(b-c)(b-d)(b-e)(b-f)\;\textrm{etc.}&=&B,\\
(c-a)(c-b)(c-d)(c-e)(c-f)\;\textrm{etc.}&=&C,\\
(d-a)(d-b)(d-c)(d-e)(d-f)\;\textrm{etc.}&=&D,\\
(e-a)(e-b)(e-c)(e-d)(e-f)\;\textrm{etc.}&=&E\\
\textrm{etc.}&&
\end{eqnarray*}
{\em each of which is made from $m-1$ factors, then not only will it
be as before}
\[
\frac{1}{A}+\frac{1}{B}+\frac{1}{C}+\frac{1}{D}+\frac{1}{E}+\textrm{etc.}=0
\]
{\em but even in this general way}
\[
\frac{a^n}{A}+\frac{b^n}{B}+\frac{c^n}{C}+\frac{d^n}{D}+\frac{e^n}{E}+\textrm{etc.}=0,
\]
{\em providing that the exponent $n$ is a positive integral number less than $m-1$.}

Thus in the above  example where the given numbers
are $3,8,12,15,17,18$, not only is it there
\[
\frac{1}{113400}-\frac{1}{12600}+\frac{1}{3240}-\frac{1}{1512}+\frac{1}{1260}
-\frac{1}{2700}=0,
\]
but one even obtains the truth of the following fractions
\begin{eqnarray*}
\frac{3}{113400}-\frac{8}{12600}+\frac{12}{3240}-\frac{15}{1512}
+\frac{17}{1260}-\frac{18}{2700}&=&0,\\
\frac{3^2}{113400}-\frac{8^2}{12600}+\frac{12^2}{3240}-\frac{15^2}{1512}
+\frac{17^2}{1260}-\frac{18^2}{2700}&=&0,\\
\frac{3^3}{113400}-\frac{8^3}{12600}+\frac{12^3}{3240}-\frac{15^3}{1512}
+\frac{17^3}{1260}-\frac{18^3}{2700}&=&0,\\
\frac{3^4}{113400}-\frac{8^4}{12600}+\frac{12^4}{3240}-\frac{15^4}{1512}
+\frac{17^4}{1260}-\frac{18^4}{2700}&=&0;
\end{eqnarray*}
but one is not permitted to continue this to higher powers, since each denominator
is made from five factors.

\begin{center}
{\Large Demonstration of the theorem\footnote{Translator: cf. Gauss, {\em Carl Friedrich Gauss Werke, Band III},
pp. 265--268.}}
\end{center}

I found this theorem from the consideration of the formula\footnote{Translator:
The denominator here has only $m-1$ factors: it is
$(x-a)(x-b)(x-c)(x-d)\cdots(x-\upsilon)$, where
$\upsilon$ is the second last of the numbers. For the moment $x$
is a variable, but later we will take $x$ to be the last
of the numbers.}
\[
\frac{x^n}{(x-a)(x-b)(x-c)(x-d)\;\textrm{etc.}},
\]
which, whenever the exponent $n$ is a positive integral
numbers less than the number of factors in the denominator,
is such that it can always be resolved into simple fractions
as such\footnote{Translator: Partial fractions: see Euler's {\em Introductio
in analysin infinitorum}, vol. I, \S 46.}
\[
\frac{A'}{x-a}+\frac{B'}{x-b}+\frac{C'}{x-c}+\frac{D'}{x-d}+\textrm{etc.},
\] 
whose denominators are the very factors of the denominator, and whose
numerators are constant quantities, not depending on $x$,
each of which can be defined in the following way.
Since the given form is equal to these simple fractions,
by multiplying by $x-a$ we will have
\[
\frac{x^n}{(x-b)(x-c)(x-d)\;\textrm{etc.}}=A'+\frac{B'(x-a)}{x-b}
+\frac{C'(x-a)}{x-c}+\frac{D'(x-a)}{x-d}+\textrm{etc.}
\]
This equality holds for any value taken for $x$, since the letters
$A',B',C',D'$, etc. do not depend on $x$. Therefore this equation will be
true if one takes $x=a$, whence
\[
\frac{a^n}{(a-b)(a-c)(a-d)\;\textrm{etc.}}=A',
\]
and thus the value of $A'$ is known. One similarly sees that
\[
B'=\frac{b^n}{(b-a)(b-c)(b-d)\;\textrm{etc.}},\quad
C'=\frac{c^n}{(c-a)(c-b)(c-d)\;\textrm{etc.}},
\]
and so on for the others. Then, moving the simple fractions to the one
side,
\[
\frac{x^n}{(x-a)(x-b)(x-c)(x-d)\;\textrm{etc}}
+\frac{A'}{a-x}+\frac{B'}{b-x}+\frac{C'}{c-x}+\frac{D'}{d-x}
+\textrm{etc.}=0,
\]
and in any case we will have, viewing the number $x$ as the last of the
numbers $a,b,c,d,\ldots,x$,
\[
\begin{split}
&\frac{a^n}{(a-b)(a-c)(a-d)\cdots(a-x)}+\frac{b^n}{(b-a)(b-c)(b-d)\cdots(b-x)}\\
&+\frac{c^n}{(c-a)(c-b)(c-d)\cdots(c-x)}+\ldots
+\frac{x^n}{(x-a)(x-b)(x-c)\cdots(x-\upsilon)}=0,
\end{split}
\]
with $\upsilon$ denoting the second last of the numbers.

This is the demonstration of the proposed theorem, which is thus not
quite obvious,
so that it seems that this truth should be counted among the common ones,
whose rule is easily perceived, unless perhaps another easier
demonstration could be found;\footnote{Translator: Then
not only would the statement of the proof be easily perceived,
but also its proof would be easily perceived?}
but the nature of the rule hardly lets one hope for this, because
this theorem is not true unless the exponent $n$ is a positive integral
number less than the number of factors in each of the denominators.

Then, since taking a greater number for $n$, the sum of these fractions
no longer vanishes, from the same source from which we drew this
theorem, for each case we will be able to assign the value
of the sum, namely by taking the number of factors to be $=m-1$
and  
hence the number of all the given numbers $a,b,c,d,\ldots,x$ to be
$=m$. If $n=m-1$, or $n=m$, or $n>m$, the fraction  
\[
\frac{x^n}{(x-a)(x-b)(x-c)(x-d)\;\textrm{etc.}}
\]
used in the demonstration should be seen as improper,
and contains as it were an integral parts,\footnote{Translator:
cf. vol. I, \S 46 of the {\em Introductio}.
If $\deg M \geq \deg N$ then
there is a polynomial part in the partial
fraction decomposition of $\frac{M}{N}$. In \S 38 Euler defines
an ``improper'' fraction of polynomials.}
and the
sum of the fractions will be equal to that very part.

Thus in the case in which $n=m-1$, the integral part is unity, and so
too the sum of the fractions $=1$. Hence in the example treated above,
where the signs are changed according to the demonstration, it will be 
\[
\frac{18^5}{2700}-\frac{17^5}{1260}+\frac{15^5}{1512}
-\frac{12^5}{3240}+\frac{8^5}{12600}-\frac{3^5}{113400}=1.
\]

But if $n=m$, the integral part arising from the fraction is
\[
x+a+b+c+d+\textrm{etc.},
\]
that is, the sum of all the given numbers. 
Therefore since in the above example the sum of all
the given numbers is $=73$, it will be
\[
\frac{18^6}{2700}-\frac{17^6}{1260}+\frac{15^6}{1512}-\frac{12^6}{3240}
+\frac{8^6}{12600}-\frac{3^6}{113400}=73.
\]

One can easily see from here how the further sums are to be found.
Namely, first the sum of all the given numbers $a,b,c,d,\ldots,x$ is taken,
which one lets $=P$, then the sum of the products from two, which one lets
$=Q$, next the sum of the products from three, which one lets $=R$,
likewise from four $=S$, from five $=T$, and so on. Now with this done,
one forms the series\footnote{Translator: Rather one forms
the generating function $1+\mathfrak{A}z+\mathfrak{B}z^2+
\mathfrak{C}z^3+\mathfrak{D}z^4+$ etc.}
\[
1+\mathfrak{A}+\mathfrak{B}+\mathfrak{C}+\mathfrak{D}+\textrm{etc.}
\]
such that
\[
\begin{split}
&\mathfrak{A}=P, \quad \mathfrak{B}=\mathfrak{A}P-Q,\quad
\mathfrak{C}=\mathfrak{B}P-\mathfrak{A}Q+R,\\
&\mathfrak{D}=\mathfrak{C}P-\mathfrak{B}Q+\mathfrak{A}R-S \quad
\textrm{etc.}
\end{split}
\]
and then
\[
\begin{array}{l|l}
\textrm{case}&\textrm{sum of the fractions}\\
n=m-1&1,\\
n=m&\mathfrak{A}=P,\\
n=m+1&\mathfrak{B}=P^2-Q,\\
n=m+2&\mathfrak{C}=P^3-2PQ+R,\\
n=m+3&\mathfrak{D}=P^4-3P^2Q+2PR+Q^2-S,\\
n=m+4&\mathfrak{E}=P^5-4P^3Q+3P^2R+3PQ^2-2PS-2QR+T\\
&\textrm{etc.}
\end{array}
\]

Or, if one puts the sum of the numbers $=\mathfrak{P}$, the sum 
of their squares $=\mathfrak{Q}$,
the sum of their cubes $=\mathfrak{R}$, the sum of their fourth powers
$=\mathfrak{S}$, of their fifth powers $=\mathfrak{T}$, etc., it will
be such that\footnote{Translator: The {\em Opera omnia} refers to Euler's E153,
{\em Demonstratio gemina theorematis Neutoniani\ldots}, in which Euler
proves Newton's identities. Newton's identities relate the coefficients
of a polynomial with the sums of the powers of the roots of the polynomial.}
\[
\begin{split}
&\mathfrak{A}=\mathfrak{P},\quad \mathfrak{B}=\frac{1}{2}\mathfrak{P}^2+\frac{1}{2}\mathfrak{Q},
\quad \mathfrak{C}=\frac{1}{6}\mathfrak{P}^3+\frac{1}{2}\mathfrak{P}\mathfrak{Q}
+\frac{1}{3}\mathfrak{R},\\
&\mathfrak{D}=\frac{1}{24}\mathfrak{P}^4+\frac{1}{4}\mathfrak{P}^2\mathfrak{Q}
+\frac{1}{8}\mathfrak{Q}^2+\frac{1}{3}\mathfrak{P}\mathfrak{R}
+\frac{1}{4}\mathfrak{S} \quad \textrm{etc.},
\end{split}
\]
which values proceed according to the law
\begin{eqnarray*}
\mathfrak{A}&=&\mathfrak{P},\\
\mathfrak{B}&=&\frac{1}{2}(\mathfrak{P}\mathfrak{A}+\mathfrak{Q}),\\
\mathfrak{C}&=&\frac{1}{3}(\mathfrak{P}\mathfrak{B}+\mathfrak{Q}\mathfrak{A}
+\mathfrak{R}),\\
\mathfrak{D}&=&\frac{1}{4}(\mathfrak{P}\mathfrak{C}+\mathfrak{Q}\mathfrak{B}
+\mathfrak{R}\mathfrak{A}+\mathfrak{S})\\
&&\textrm{etc.}
\end{eqnarray*}

With the truth of our theorem established, I judge that it will not be
beside the point if I should carefully investigate the nature of the
formulas on which the theorem turns. Therefore, if numbers $a,b,c,d$, etc.
are given, for each of them one searches for what the character will be of
the formula $(a-b)(a-c)(a-d)$ etc., which is produced from the product
of the differences
of it with all the others.
Then let the number of the given numbers be $=n$,\footnote{Translator: This
was denoted earlier by $m$.} and assuming $z$ to be a variable quantity,
I form this product from it
\[
(z-a)(z-b)(z-c)(z-d)(z-e)\;\textrm{etc.},
\]
which by multiplication yields the expansion
\[
z^n-Pz^{n-1}+Qz^{n-2}-Rz^{n-3}+Sz^{n-4}-\textrm{etc.}
\]
Therefore by dividing by $z-a$ we will have
\[
(z-b)(z-c)(z-d)\;\textrm{etc.}=\frac{z^n-Pz^{n-1}+Qz^{n-2}-Rz^{n-3}+\textrm{etc.}}{z-a}.
\]
If we now put here $z=a$, the very form $(a-b)(a-c)(a-d)\;\textrm{etc.}$ will
arise which I indicated above with the letter $A$. Then indeed, for the
other side, both the numerator and the denominator go to zero, and therefore
its value will be
\[
na^{n-1}-(n-1)Pa^{n-2}+(n-2)Qa^{n-3}-(n-3)Ra^{n-4}+\textrm{etc.},
\]
which, since
it is the case that
\[
a^n-Pa^{n-1}+Qa^{n-2}-Ra^{n-3}+Sa^{n-4}-\textrm{etc.}=0,
\]

\ldots\footnote{Translator: The paper stops here.}

\end{document}